\documentclass[11pt]{amsart}
\usepackage{amsmath,amssymb,latexsym,amsfonts,amscd}
\def\m{\mathfrak m}
\def\inc{\subset}
\newtheorem{theorem}{Theorem}[section]

\newtheorem{lemma}[theorem]{Lemma}
\newtheorem{corollary}[theorem]{Corollary}

\begin{document}

\title{Absolute integral closure in positive characteristic }
\author{Craig Huneke}
\address{Department of Mathematics, University of Kansas, Lawrence, KS 66045}
\email{huneke@math.ku.edu}
\author{Gennady Lyubeznik}
\thanks{NSF support for both authors is gratefully acknowledged. The first author was supported in part
by grant DMS-0244405, and the second by grant DMS-0202176.}
\address{Department of Mathematics, University of Minnesota, Minneapolis,
MN 55455}
\email{gennady@math.umn.edu}
\begin{abstract} Let $R$ be a local Noetherian domain of positive characteristic.
 A theorem of Hochster and Huneke (1992) states that if $R$ is excellent, then
the absolute integral closure of $R$ 
is a big Cohen-Macaulay algebra. We prove that if $R$ is 
the homomorphic image of a Gorenstein local ring, then all the local cohomology (below the dimension) of such a ring
maps to zero in a finite extension of the ring. There results an extension of the original result
of Hochster and Huneke to the case in which $R$ is a homomorphic image of a Gorenstein
local ring, and  a considerably 
simpler proof of this result in
the cases where the assumptions overlap, e.g., for complete Noetherian local domains.\end{abstract}
\date{March 31, 2006}



\maketitle

\section{Introduction}

Let $R$ be a commutative Noetherian domain with fraction field $K$. The \it absolute integral closure \rm
of $R$, denoted $R^+$, is the integral closure of $R$ in a fixed algebraic closure  $\overline{K}$ of $K$.
This ring was studied by Artin in \cite{Ar} where among other results he proved that
in the case $R$ is Henselian and local, the sum of two primes ideals of $R^+$ remains prime.

In \cite{HH}, Hochster and Huneke proved that if $(R,\m)$ is an excellent local Noetherian domain of
positive characteristic $p > 0$,
then $R^+$ is a big Cohen-Macaulay algebra, i.e., every system of parameters in $R$ is a regular sequence
on $R^{+}$. The corresponding statement in equicharacteristic $0$ is false if the dimension
is at least three. 
Smith \cite{Sm} further proved that the tight closure of an ideal $I$ generated by  parameters
is exactly the extension and  contraction of
$I$ to $R^+$: $I^* = IR^+\cap R$. It is an open question whether the latter equality is true
for every ideal $I$ in an excellent Noetherian local domain of positive characteristic. 
See \cite{A},  \cite{AH}, and \cite{Si} for additional work concerning $R^+$.

Throughout this paper, $R$ is a commutative Noetherian domain of characteristic $p > 0$ with fraction field
$K$, $\overline{K}$ is a
fixed algebraic closure of $K$, and $R^{+}$ is the integral closure of $R$ in $\overline{K}$.
The theorem of Hochster and Huneke in \cite{HH} implies the following:

\smallskip

{\it Let $(R, \m)$ be an excellent local commutative Noetherian domain of characteristic $p>0$.
Then the  natural  homomorphism $ H^i_{\m}(R)\rightarrow H^i_{\m}(R^{+})$ is the
zero map for every $i < \dim R$.}

\smallskip 

In fact, as we observe in this paper (see Corollary~\ref{corcm}), the statement above is basically equivalent to the statement that
$R^+$ is a big Cohen-Macaulay algebra for $R$.
 
The main result of this paper, Theorem~\ref{module} below,  states that if $R$ is a homomorphic
image of a Gorenstein local ring (though not necessarily excellent), then in fact one can find
a \emph{finite} extension ring $S$,
$R\inc S\inc R^+$, such that the map from $H^i_{\m}(R)\rightarrow H^i_{\m}(S)$ is zero for all
$i < \dim R$. Our proof is independent of the results of \cite{HH} and in particular gives a
considerably simpler proof of the main result of that paper, with a stronger conclusion, when
the assumptions overlap. For example, if $R$ is complete, then it is both excellent and
a homomorphic image of a Gorenstein ring.
Our proof was in part inspired by the work of Hartshorne and Speiser in \cite{HS} and Lyubeznik in \cite{Ly},
concerning the structure of local cohomology modules in positive characteristic.

\bigskip

\section{Main Result}

\medskip

Let $R$ be a commutative ring containing a field of characteristic $p>0$, let $I\subset R$
be an ideal, and let $R'$ be an $R$-algebra. The Frobenius ring homomorphism
$f:R'\stackrel{r\mapsto r^p}{\to}R'$ induces a map $f_*:H^i_I(R')\to H^i_I(R')$ on all
local cohomology modules of $R'$ called the action of the Frobenius on $H^i_I(R')$.
For an element $\alpha\in H^i_I(R')$ we denote $f_*(\alpha)$ by $\alpha^p$.

\medskip

We recall that for a Gorenstein local ring $A$ of dimension $n$, local duality says that
there is an isomorphism of functors $D({\rm Ext}^{n-i}_A(-, A))\cong H^i_{\mathfrak m}(-)$
on the category of finite $A$-modules, where $D={\rm Hom}_A(-,E)$ is the Matlis duality
functor (here $E$ is the injective hull of the residue field of $A$ in the category of
$A$-modules) \cite[11.2.6]{BS}.

\begin{theorem}\label{module}
Let $R$ be a commutative Noetherian local domain containing a field of characteristic $p>0$, let
$K$ be the fraction field of $R$ and let $\overline K$ be the algebraic closure of $K$. Assume $R$ is a surjective
image of a Gorenstein local ring $A$. Let $\mathfrak m$ be the maximal ideal of $R$. Let $R'$ be an $R$-subalgebra
of $\overline K$ (i.e. $R\subset R'\subset \overline K$) that is a finite $R$-module. Let $i< \dim R$
be a non-negative integer. There is an $R'$-subalgebra $R''$ of $\overline K$ (i.e. $R'\subset R''\subset \overline K$)
that is finite as an $R$-module and such that the natural map $H^i_{\mathfrak m}(R')\to H^i_{\mathfrak m}(R'')$
is the zero map.
\end{theorem}

\emph{Proof.} Let $n=\dim A$ and let $N={\rm Ext}^{n-i}_A(R',A)$. Since $R'$ is a finite $R$-module,
so is $N$. 

Let $d=\dim R$. We use induction on $d$. For $d=0$ there is nothing to prove, so we assume that
$d>0$ and the theorem proven for all smaller dimensions. Let $P\subset R$ be a non-maximal prime ideal.
We claim there exists an $R'$-subalgebra $R^P$ of $\overline K$ (i.e. $R'\subset R^P\subset \overline K$) such that
$R^P$ is a finite $R$-module and for every $R^P$-subalgebra $R^*$ of $\overline K$ (i.e. $R^P\subset R^*\subset \overline K$)
such that $R^*$ is a finite $R$-module, the image $\mathcal I\subset N$ of the natural map
${\rm Ext}^{n-i}_A(R^*,A)\to N$ induced by the natural inclusion $R'\to R^*$ vanishes after localization at $P$,
i.e. $\mathcal I_P=0$. Indeed, let $d_P=\dim R/P$. Since $P$ is different from the maximal ideal,
$d_P>0$. As $R$ is a surjective image of a Gorenstein local ring, it is catenary, hence the dimension
of $R_{P}$ equals $d-d_P$, and $i<d$ implies $i-d_P<d-d_P=\dim R_{P}$. By the induction hypothesis
applied to the local ring $R_{P}$ and the $R_{P}$-algebra $R'_{P}$, which is finite as an $R_{P}$-module,
there is an $R'_{P}$-subalgebra $\tilde R$ of $\overline K$, which is finite as an $R_{P}$-module, such that
the natural map $H^{i-d_P}_{P}(R'_P)\to H^{i-d_P}_{P}(\tilde R)$ is the zero map.
Let $\tilde R=R'_{P}[z_{1},z_{2},\dots,z_{t}]$, where $z_{1},z_{2},\dots,z_{t}\in \overline K$ are
integral over $R_{P}$. Multiplying, if necessary, each $z_{j}$ by some element of $R\setminus P$,
we can assume that each $z_{j}$ is integral over $R$. We set $ R^P=R'[z_{1},z_{2},\dots,z_{t}]$.
Clearly, $R^P$ is an $R'$-subalgebra of $\overline K$ that is finite as $R$-module. 

Now let $R^*$ be both
an $R^P$-subalgebra of $\overline K$ (i.e. $R^P\subset R^*\subset \overline K$) and a finite $R$-module.
The natural inclusions $R'\to R^P\to R^*$ induce natural maps 
${\rm Ext}^{n-i}_A(R^*,A)\to {\rm Ext}^{n-i}_A(R^P,A)\to N$. This implies that
$\mathcal I\subset \mathcal J$, where $\mathcal J$ is the image of the natural map
$\phi:{\rm Ext}^{n-i}_A(R^P,A)\to N$. Hence it is enough to prove that $\mathcal J_{P}=0$.
Localizing this map at $P$ we conclude that $J_P$ is the image of the natural map
$\phi_P:{\rm Ext}^{n-i}_{A_P}(\tilde R,A_P)\to {\rm Ext}^{n-i}_{A_P}(R'_P,A_P)$ induced by the
natural inclusion $R'_P\to \tilde R$ (by a slight abuse of language we identify the prime ideal
$P$ of $R$ with its full preimage in $A$). Let $D_P(-)={\rm Hom}_{A_P}(-, E_P)$ be the Matlis duality
functor in the category of $R_P$-modules, where $E_P$ is the injective hull of the residue field of
$R_P$ in the category of $R_P$-modules. Local duality implies that $D_P(\phi_P)$ is the natural map
$H^{i-d_P}_{P}(R'_P)\to H^{i-d_P}_{P}(\tilde R)$ which is the zero map by construction
(note that $i-d_P=\dim A_P-(n-i)$). Since $\phi_P$ is a map between finite $R_P$-modules
and $D_P(\phi_P)=0$, it follows that $\phi_P=0$. This proves the claim.

Since $N$ is a finite $R$-module, the set of the associated primes of $N$ is finite.
Let $P_1,\dots,P_s$ be the associated primes of $N$ different from $\mathfrak m$. For each $j$
let $R^{P_j}$ be an $R'$-subalgebra of $\overline K$ corresponding to $P_j$, whose existence is guaranteed
by the above claim. Let $\overline R'=R'[R^{P_1},\dots, R^{P_s}]$ be the compositum of all the $R^{P_j}$, $1\leq j\leq s$.
Clearly, $\overline R'$ is an $R'$-subalgebra of $\overline K$ (i.e. $R'\subset \overline R'\subset \overline K$).  Since each $R^{P_j}$
is a finite $R$-module, so is $\overline R'$. Clearly, $\overline R'$  contains every $R^{P_j}$. Hence the above
claim implies that $\mathcal I_{P_j}=0$ for every $j$, where $\mathcal I\subset N$ is the image of the
natural map ${\rm Ext}^{n-i}_A(\overline R',A)\to N$ induced by the natural inclusion $R'\to\overline R'$. It follows
that  not a single $P_j$ is an associated prime of $\mathcal I$. But $\mathcal I$ is a submodule of $N$,
and therefore every associated prime of $\mathcal I$ is an associated prime of $N$.
Since $P_1,\dots, P_s$ are all the associated primes of $N$ different from $\mathfrak m$, we conclude that if
$\mathcal I\ne 0$, then $\mathfrak m$ is the only associated prime of $\mathcal I$. Since $\mathcal I$,
being a submodule of a finite $R$-module $N$, is finite, and since $\mathfrak m$ is the only associated
prime of $\mathcal I$, we conclude that $\mathcal I$ is an $R$-module of finite length.

Writing the natural map ${\rm Ext}^{n-i}_A(\overline R',A)\to N$ as the composition of two maps
${\rm Ext}^{n-i}_A(\overline R',A)\to \mathcal I\to N$, the first of which is surjective and the second injective,
and applying the Matlis duality functor $D$, we get that the natural map
$\varphi:H^i_{\mathfrak m}(R')\to H^i_{\mathfrak m}(\overline R')$ induced by the inclusion $R'\to \overline R'$ is the
composition of two maps $H^i_{\mathfrak m}(R')\to D(\mathcal I)\to H^i_{\mathfrak m}(\overline R')$, the first
of which is surjective and the second injective. This shows that the image of $\varphi$ is isomorphic to
$D(\mathcal I)$ which is an $R$-module of finite length since so is $\mathcal I$. In particular, the image
of $\varphi$ is a finitely generated $R$-module. Let $\alpha_1,\dots, \alpha_s\in
H^i_{\mathfrak m}(\overline R')$ generate Im$\varphi$.

The natural inclusion $R'\to \overline R'$ is compatible with the Frobenius homomorphism, i.e. with
the raising to the $p$th power on $R'$ and $\overline R'$. This implies that $\varphi$ is compatible
with the action of the Frobenius $f_*$ on $H^i_{\mathfrak m}(R')$ and $H^i_{\mathfrak m}(\overline R')$,
i.e. $\varphi(f_*(\alpha))=f_*(\varphi(\alpha))$ for every $\alpha\in H^i_{\mathfrak m}(R')$,
which, in turn, implies that Im$\varphi$ is an $f_*$-stable $R$-submodule of $H^i_{\mathfrak m}(\overline R')$,
i.e. $f_*(\alpha)\in {\rm Im}\varphi$ for every $\alpha\in {\rm Im}\varphi$. We finish the proof
by applying the following lemma to each element of  a finite generating set $\alpha_1,...,
\alpha_s$ of Im$\varphi$.
Applying Lemma~\ref{element} below we obtain a
$\overline R'$-subalgebra $R'_j$ of $\overline K$ (i.e. $\overline R'\subset R'_j\subset \overline K$) such that
$R'_j$ is a finite $R$-module and the natural map $H^i_{\mathfrak m}(\overline R')\to H^i_{\mathfrak m}(R'_j)$
sends $\alpha_j$ to zero. Let $R''=R'[R_1,\dots, R_s]$ be the compositum of all the $R'_j$.
Then $R''$ is an $R'$-subalgebra of $\overline K$ and is a finite $R$-module since so is each $R'_j$.
The natural map $H^i_{\mathfrak m}(\overline R')\to H^i_{\mathfrak m}(R'')$ sends every $\alpha_j$ to zero,
hence it sends the entire Im$\varphi$ to zero. Thus the natural map
$H^i_{\mathfrak m}(R')\to H^i_{\mathfrak m}(R'')$ is zero. \qed

\medskip

To finish the proof we prove the following lemma, which is closely related to the
``equational lemma" in \cite{HH} and its modification in \cite{Sm}, (5.3).

\begin{lemma} \label{element} Let $R$ be a commutative Noetherian domain containing a field of
characteristic $p>0$, let $K$ be the fraction field of $R$ and let $\overline K$ be the algebraic closure
of $K$. Let $I$ be an ideal of $R$ and let $\alpha\in H^i_I(R)$ be an element such that the elements
$\alpha, \alpha^p,\alpha^{p^2},\dots,\alpha^{p^t},\dots$ belong to a finitely generated $R$-submodule of $H^i_I(R)$.
There exists an $R$-subalgebra $R'$ of $\overline K$ (i.e. $R\subset R'\subset \overline K$) that is finite
as an $R$-module and such that the natural map $H^i_I(R)\to H^i_I(R')$ induced by the natural
inclusion $R\to R'$ sends $\alpha$ to 0.
\end{lemma}
                                                                                          
\emph{Proof.} Let $A_t=\sum_{i=1}^{i=t}R\alpha^{p^i}$ be the $R$-submodule of $H^i_I(R)$ generated
by $\alpha,\alpha^p,\dots,\alpha^{p^t}$. The ascending chain $A_1\subset A_2\subset A_3\subset\dots$
stabilizes because $R$ is Noetherian and all $A_t$ sit inside a single finitely generated $R$-submodule
of $H^i_I(R)$. Hence $A_s=A_{s-1}$ for some $s$, i.e. $\alpha^{p^s}\in A_{s-1}$. Thus there exists an
equation $\alpha^{p^s}=r_1\alpha^{p^{s-1}}+r_2\alpha^{p^{s-2}}+\dots+r_{s-1}\alpha$ with $r_i\in R$
for all $i$. Let $T$ be a variable and let $g(T)=T^{p^s}-r_1T^{p^{s-1}}-r_2^{p^{s-2}}-\dots-r_{s-1}T$.
Clearly, $g(T)$ is a monic polynomial in $T$ with coefficients in $R$ and $g(\alpha)=0$.
                                                                                          
Let $x_1,\dots, x_d\in R$ generate the ideal $I$. If $M$ is an $R$-module, the \v Cech complex 
$C^{\bullet}(M)$ of $M$ with respect to the generators $x_1,\dots, x_d\in R$ is
$$0\to C^0(M)\to\dots \to C^{i-1}(M)\stackrel{d_{i-1}}{\to} C^i(M)\stackrel{d_i}{\to}
C^{i+1}(M)\to\dots\to C^d(M)\to 0$$
where $C^0(M)=M$ and $C^i(M)=\oplus_{1\leq j_1<\dots<j_{i}\leq d}R_{x_{j_1}\cdots x_{j_{i}}}$,
and $H^i_I(M)$ is the $i$th cohomology module of $C^{\bullet}(M)$ \cite[5.1.19]{BS}.
                                                                                          
Let $\tilde \alpha\in C^i(R)$ be a cycle (i.e. $d_i(\tilde \alpha)=0$) that represents $\alpha$.
The equality $g(\alpha)=0$ means that $g(\tilde \alpha)=d_{i-1}(\beta)$ for some $\beta\in C^{i-1}(R)$.
Since $C^{i-1}(R)=\oplus_{1\leq j_1<\dots<j_{i-1}\leq d}R_{x_{j_1}\cdots x_{j_{i-1}}}$, we may write
$\beta=(\frac{r_{j_1,\dots,j_{i-1}}}{x_{j_1}^{e_{1}}\cdots x_{j_{i-1}}^{e_{i-1}}})$
where $r_{j_1,\dots,j_{i-1}}\in R$, the integers $e_1,\dots, e_{i-1}$ are non-negative, and
$\frac{r_{j_1,\dots,j_{i-1}}}{x_{j_1}^{e_{1}}\cdots x_{j_{i-1}}^{e_{i-1}}}\in R_{x_{j_1}\cdots x_{j_{i-1}}}$.
                                                                                          
Consider the equation $g(\frac{Z_{j_1,\dots,j_{i-1}}}{x_{j_1}^{e_{1}}\cdots x_{j_{i-1}}^{e_{i-1}}})
-\frac{r_{j_1,\dots,j_{i-1}}}{x_{j_1}^{e_{1}}\cdots x_{j_{i-1}}^{e_{i-1}}}=0$ where
$Z_{j_1,\dots,j_{i-1}}$ is a variable. Multiplying this equation by 
$(x_{j_1}^{e_{1}}\cdots x_{j_{i-1}}^{e_{i-1}})^{p^s}$ produces a monic polynomial equation in
$Z_{j_1,\dots,j_{i-1}}$ with coefficients in $R$. Let $z_{j_1,\dots,j_{i-1}}\in \overline K$ be a root
of this equation and let $R''$ be the $R$-subalgebra of $\overline K$ generated by all the $z_{j_1,\dots,j_{i-1}}$s,
i.e. by the set $\{z_{j_1,\dots,j_{i-1}}|1\leq j_1<\dots<j_{i-1}\leq d\}$. Since each $z_{j_1,\dots,j_{i-1}}$
is integral over $R$ and there are finitely many $z_{j_1,\dots,j_{i-1}}$s, the $R$-algebra $R''$ is finite
as an $R$-module.

Let $\tilde{\tilde\alpha} =
(\frac{z_{j_1,\dots,j_{i-1}}}{x_{j_1}^{e_{1}}
\cdots x_{j_{i-1}}^{e_{i-1}}})\in C^{i-1}(R'')$. The natural inclusion $R\to R''$ makes $C^{\bullet}(R)$
into a subcomplex of $C^{\bullet}(R'')$ in a natural way, and we identify $\tilde \alpha\in C^i(R)$
and $\beta\in C^{i-1}(R)$ with their natural images in $C^i(R'')$ and $C^{i-1}(R'')$ respectively.
With this identification, $\tilde \alpha\in C^i(R'')$ is a cycle representing the image of $\alpha$
under the natural map $H^i_I(R)\to H^i_I(R'')$, and so is $\overline \alpha=\tilde\alpha-d_{i-1}(\tilde{\tilde\alpha})
\in C^i(R'')$. Since $g(\tilde{\tilde\alpha})=\beta$ and $g(\tilde \alpha)=d_{i-1}(\beta)$, we conclude that
$g(\overline \alpha)=0$. Let $\overline\alpha=(\rho_{j_1,\dots,j_i})$ where $\rho_{j_1,\dots,j_i}
\in R''_{x_{j_1}\cdots x_{j_i}}$. Each individual $\rho_{j_1,\dots,j_i}$ satisfies the equation
$g(\rho_{j_1,\dots,j_i})=0$. Since $g(T)$ is a monic polynomial in $T$ with coefficients in $R$,
each $\rho_{j_1,\dots,j_i}$ is an element of the fraction field of $R''$ that is integral over $R$.
Let $R'$ be obtained from $R''$ by adjoining all the $\rho_{j_1,\dots,j_i}$.
                                                                                          
Each $\rho_{j_1,\dots,j_i}\in R'$, so the image of $\alpha$ in $H^i_I(R')$ is represented by the cycle
$\overline \alpha=(\rho_{j_1,\dots,j_i})\in C^i(R')$ which has all its components $\rho_{j_1,\dots,j_i}$ in $R'$. 
Each $R'_{x_{j_1}\cdots x_{j_i}}$ contains a natural copy of $R'$, namely, the one generated by the element
$1\in R'_{x_{j_1}\cdots x_{j_i}}$. There is a subcomplex of $C^{\bullet}(R')$ that in each degree is the
direct sum of all such copies of $R'$. This subcomplex is exact because its cohomology groups are the
cohomology groups of $R'$ with respect to the unit ideal. Since $\overline \alpha$ is a cycle and belongs to
this exact subcomplex, it is a boundary, hence it represents the zero element in $H^i_I(R')$. \qed

\begin{corollary}\label{corcm}
Let $R$ be a commutative Noetherian local domain containing a field of characteristic $p>0$.
Assume that $R$ is a surjective image of a Gorenstein local ring. Then the following hold:

(a) $H^i_{\mathfrak m}(R^+)=0$ for all $i<\dim R$, where $\mathfrak m$ is the maximal ideal of $R$.

(b) Every system of parameters of $R$ is a regular sequence on $R^+$.
\end{corollary}

\emph{Proof.} (a) $R^+$ is the direct limit of the finitely generated $R$-subalgebras $R'$, hence $H^i_{\mathfrak m}(R^+)=\varinjlim H^i_{\mathfrak m}(R')$. But Theorem \ref{module} implies that for each $R'$ there is $R''$ such that the map $H^i_{\mathfrak m}(R')\to H^i_{\mathfrak m}(R'')$ in the inductive system is zero. Hence the limit is zero.

(b) Let $x_1,..., x_d$ be a system of parameters of $R$. We prove that $x_1,...,x_j$ is a regular
sequence on $R^+$ by induction on $j$. The case $j=1$ is clear, since $R^+$ is a domain. 
Assume that $j>1$ and $x_1,\dots, x_{j-1}$ is a regular sequence on $R^+$. Set $I_t = (x_1,...,x_t)$.
The fact that $H^i_{\mathfrak m}(R^+)=0$ for all $i<d$ and the short exact sequences 
$$0\to R^+/I_{t-1}R^+\stackrel{\rm mult\ by\ x_t}{\longrightarrow} R^+/I_{t-1}R^+\to R^+/I_{t}R^+\to 0$$ for
$t\leq j-1$ imply by induction on $t$ that $H^q_{\mathfrak m}(R^+/(x_1,...,x_{t})R^+)=0$ for $q<d-t$.
In particular, $H^0_{\mathfrak m}(R^+/(x_1,...,x_{j-1})R^+)=0$ since $0<d-(j-1)$.
Hence $\mathfrak m$ is not an associated prime of $R^+/(x_1,...,x_{j-1})R^+$. This implies that the
only associated primes of $R^+/(x_1,...,x_{j-1})R^+$ are the minimal primes of $R/(x_1,...,x_{j-1})R$.
Indeed, if there is an embedded associated prime, say $P$, then $P$ is the maximal ideal of the ring
$R_P$ whose dimension is bigger than $j-1$ and $P$ is an associated prime of 
$(R^+/(x_1,...,x_{j-1})R^+)_P=(R_P)^+/(x_1,...,x_{j-1})(R_P)^+$ which is impossible by the above. 
Hence every element of $\mathfrak m$ not in any minimal prime of $R/(x_1,...,x_{j-1})R$, for example,
$x_j$, is a regular element on $R^+/(x_1,...,x_{j-1})R^+$.\qed

\end{document}